\documentclass[a4paper]{amsart}

\usepackage{amssymb}
\usepackage[colorlinks=true, linkcolor=blue]{hyperref}
\usepackage[all]{xypic}

\swapnumbers

\makeatletter
\def\@secnumfont{\bfseries}
\makeatletter

\theoremstyle{thm}
\newtheorem{thm}[subsection]{Theorem}
\newtheorem{prop}[subsection]{Proposition}
\newtheorem*{prop*}{Proposition}
\newtheorem*{thm*}{Theorem}

\newtheorem*{lem*}{Lemma}

\newtheorem*{kor*}{Corollary}

\theoremstyle{definition} 

\def\frak{\mathfrak}
\def\Bbb{\mathbb}
\def\Cal{\mathcal}

\newcommand{\id}{\operatorname{id}}

\newcommand{\End}{\operatorname{End}}

\newcommand{\gr}{\operatorname{gr}}

\newcommand{\x}{\times}
\renewcommand{\o}{\circ}

\let\ccdot\cdot
\def\cdot{\hbox to 2.5pt{\hss$\ccdot$\hss}}

\newcommand{\ga}{\gamma}
\newcommand{\ka}{\kappa}

\newcommand{\om}{\omega}
\renewcommand{\phi}{\varphi}
\newcommand{\ph}{\varphi}
\newcommand{\si}{\sigma}
\renewcommand{\th}{\theta}
\newcommand{\ze}{\zeta}

\newcommand{\Ga}{\Gamma}
\newcommand{\La}{\Lambda}
\newcommand{\Ph}{\Phi}
\newcommand{\Ps}{\Psi}
\newcommand{\Rho}{{\mbox{\sf P}}}
\newcommand{\smallRho}{{\mbox{\Small\sf P}}}
\newcommand{\Om}{\Omega}

\renewcommand{\Up}{\Upsilon}
\newcommand{\cg}{{\Cal G}}
\newcommand{\ca}{{\Cal A}}
\newcommand{\cv}{{\Cal V}}
\newcommand{\bv}{{\Bbb V}}
\newcommand{\fg}{{\frak g}}
\newcommand{\fp}{{\frak p}}

\begin{document}
\title[Bundles of Weyl structures and invariant calculus]{Bundles of Weyl structures and invariant calculus for parabolic
geometries} 
\author{Andreas \v Cap and Jan Slov\'ak}

\address{A.C.: Faculty of Mathematics, University of Vienna, Oskar-Morgenstern-Platz
  1, 1090 Wien, Austria
\newline\indent J.S.: Department of Mathematics and Statistics,
Masaryk University, Kotl\'a\v rsk\'a 2a, 611~37~Brno, Czech Republic}

\email{Andreas.Cap@univie.ac.at, slovak@math.muni.cz}

\begin{abstract} For more than hundred years, various concepts were
developed to understand the fields of geometric objects and invariant differential
operators between them for conformal Riemannian and projective geometries.  More
recently, several general tools were presented for the entire class of parabolic
geometries, i.e., the Cartan geometries modelled on homogeneous spaces $G/P$ with $P$
a parabolic subgroup in a semi-simple Lie group $G$.  Similarly to conformal
Riemannian and projective structures, all these geometries determine a class of
distinguished affine connections, which carry an affine structure modelled on
differential 1-forms $\Upsilon$.  They correspond to reductions of $P$ to its
reductive Levi factor, and they are called the Weyl structures similarly to the
conformal case. The standard definition of differential invariants in
  this setting is as affine invariants of these connections, which do not depend on
  the choice within the class. In this article, we describe a universal calculus
  which provides an important first step to determine such invariants.
We present a natural procedure how to construct all affine invariants of Weyl
connections, which depend only tensorially on the deformations
$\Upsilon$.  \end{abstract}

\thanks{Dedicated to the memory of Alexandre Michailovich Vinogradov.
\\ \indent First author supported by the Austrian Science Fund (FWF): P33559-N,
second author supported by the grant GX19-28628X of GA\v CR.} 


\maketitle

Differential invariants of various geometric structures are the core ingredients for
numerous applications both in geometric analysis and mathematical physics. In
particular, the invariants of conformal Riemannian manifolds attracted a lot of
attention in the course of the last 100 years.

For smooth manifolds with an affine connection, the so called `first invariant
theorem' says that all the invariants are expressions built of the covariant
derivatives of sections of natural bundles, the curvature and the
  torsion of the connection by means of algebraic tensorial invariants,
cf. \cite{KMS} for a modern treatment. Let us call them the {\em affine differential
  invariants} on smooth manifolds.  The analogous `first invariant theorem' for
Riemannian geometries says that all differential invariants are built from affine
invariants of the canonical Levi-Civita connection via (algebraic)
  invariants of the orthogonal group.

A conformal Riemannian geometry is defined as a class of conformally equivalent
Riemannian metrics and so the above Riemannian first invariant theorem can be used to
define conformal invariants. Thus, a conformal invariant is usually understood as a
Riemannian invariant in terms of any metric from the conformal class, such that the
change of the metric does not change the invariant.  As proved by the extraordinary
effort to understand such invariants for many decades, already the first invariant
theorem is not easy in this case.

An equivalent definition of conformal structures treats them in terms of
classical G-structures as reductions of the linear frame bundle to the structure
group $G_0=CO(n)$, the group of all conformally Euclidean linear transformations in
the given dimension. Such a structure admits compatible torsion-free connections,
which are classically called \textit{Weyl connections}. This broader class of
conformal connections was exploited by H.\ Weyl which motivates their name. It turns
out that the Weyl connections form an affine space modelled on one-forms. Of course,
the Levi-Civita connections of metrics in the conformal class are Weyl connections,
they form an affine subspace modelled on exact one-forms. The study of conformal
Riemannian invariants goes back to \'E.\ Cartan, T.\ Thomas, J.A.\ Schouten, and
others (e.g., \cite{Cartan-conf, Tho1, Tho2}).  A lot of spectacular tricks to build
invariant expressions have been developed, and some of them were turned into a quite
effective calculus for conformal invariants by V.\ W\"unsch, see \cite{Wuensch}.

Motivated by the rich geometry of conformal Riemannian manifolds, the Weyl structures
and the preferred connections were introduced in the general framework of parabolic
geometries in \cite{CS} (generalizing the approach from \cite{Gau}).  In particular,
the notions of scales, closed and exact Weyl structures, and (Schouten's)
Rho-tensors were extended, and natural generalizations of classical normal
coordinates in affine geometries were discussed.

Moreover, filtered analogues of classical G-structures which are equivalent to
parabolic geometries and the general Cartan-Tanaka theory for all parabolic
geometries are explained in great detail in \cite{parabook}{, see also
  \cite{Cap:Cartan}}.  In this setting, the Weyl connections on a parabolic geometry
of type $G/P$ correspond to reductions of the parabolic structure group $P$ of the
canonical Cartan connection $\om$ on the principal bundle $\mathcal G\to M$ to its
reductive Levi factor $G_0\subset P$.

More recently, the geometry of the bundle of Weyl structures $\pi:A=\mathcal G/G_0\to
M$ was studied carefully in the joint work \cite{CapMettler} of the first author and
T.\ Mettler.  Here the canonical Cartan connection $\om$ induces a canonical affine
connection on the manifold $A$ as well as a canonical splitting $TA=L^-\oplus
L^+$. There also is a nice relation between natural bundles over $M$ associated to
$P$-representations and natural bundles over $A$. Using this, we will identify a
natural class of differential invariants of the canonical connection on $A$ which
induce affine differential invariants of the Weyl connections. The invariants
obtained in this way transform tensorially in the one-forms $\Up$ that parameterize
the Weyl connections. We shall also prove a converse to this statement, which is much
more subtle: If an affine invariant of Weyl connections, its curvature
  and its torsion transforms tensorially, then it comes from a natural affine
differential invariant on $A$ as described above.

The paper is organized as follows: After a brief review of the main tools and
concepts following \cite{parabook}, \cite{CapMettler}, we  show in
  Section 2 how affine invariants of the canonical connection $D$ on $A$ can be used
  to construct affine invariants of the Weyl connections. This is based on an
  invariant concept of jets and by construction, the resulting invariants always
  change tensorially under a change of Weyl structure. We introduce the terminology
  ``nearly invariant operators'' for this behavior. In the last section, we prove
  that all nearly invariant operators are obtained in this way. This can be viewed
  as a week version of the `first invariant theorem' (providing only an ansatz of
  possible expressions for invariant operators) and also a universal procedure
  realizing a generalization of the so-called W\"unsch calculus in conformal
  geometry, which provides universal formulae for large classes of invariant
  operators, cf.  \cite{CSS1-3}.

\section{Parabolic geometries and the bundle of Weyl structures}

We shall follow the terminology and notation of \cite{parabook} and
\cite{CapMettler}.  The complete procedure deriving the canonical (i.e., properly
normalized) Cartan connection from more elementary data on filtered manifolds
corresponding to parabolic geometries was first worked out long ago by Tanaka in
\cite{Tan}, more straightforward, simpler and more general versions appeared in
\cite{CSh}, \cite{parabook}, \cite{Cap:Cartan}. We do not go into
  details on the normalization here, but just use its implications for the curvature
  and torsion of Weyl connections. The calculus we develop actually does also work
  for non-normal Cartan connections. However in the non-normal case, in some of the
  results the notion of invariants have  to be adapted.

\subsection{The Cartan connections}\label{cartan-connection}

A Cartan geometry of type $G/P$ is an absolute parallelism on a principal fiber
bundle $\cg\to M$ with structure group $P$ encoded as a one-form
$\om\in\Om^1(\cg,\mathfrak g)$. One requires suitable equivariancy properties with
respect to the principal action of $P$, similar to those of the Maurer-Cartan form
$\om_G$ on $G\to G/P$. Thus, we may view Cartan geometries as curved deformations of
the homogeneous spaces.  More explicitly, the {\em Cartan connection}
$\om\in\Om^1(\mathcal G, \mathfrak g)$ is required to obey the following properties
\begin{itemize}
\item[(1)] $\om(\ze_X)(u) = X$ for all $X\in \fp$, $u\in\cg$ (the connection
reproduces the generators of fundamental vertical fields)
\item[(2)] $(r^g)^*\om = \operatorname{Ad}(g^{-1})\circ \om$ for any $g\in P$ (the
  connection form is equivariant with respect to the principal action)
\item[(3)] $\om_{|T_u\cg}:T_u\cg\to \fg$ is a linear isomorphism for all
$u\in \cg$ (the absolute parallelism condition).
\end{itemize}
A morphisms of Cartan geometries is a principal fiber bundle morphism
$\ph:\cg\to\cg'$ with the property $\ph^*(\om')=\om$.  {\em Parabolic geometries} are
Cartan geometries of type $(G,P)$, where $P$ is a parabolic subgroup in a semisimple
real Lie group $G$.

In the sequel, we shall consider a fixed parabolic subgroup $P\subset G$. It is well
known, that $P$ uniquely determines a grading 
\begin{equation}\label{grad}
\fg = \fg_-\oplus \fg_0\oplus \fp_+=\fg_{-k}\oplus\dots\fg_{-1}\oplus\fg_0\oplus\fg_1\oplus\dots\oplus
\fg_k
\end{equation}
on the Lie algebra $\fg$ of the group $G$, such that $\fp= \fg_0\oplus \fp_+=
\fg_0\oplus\fg_1\oplus\dots\oplus \fg_k$ is the Lie algebra of $P$. This comes with a
so called grading element $E\in\fp$ with the property $\operatorname{ad}E_{|\fg_i} =
i\cdot \operatorname{id}_{\fg_i}$ for all $i=-k,\dots,k$.

While the grading \eqref{grad} is not $P$-invariant, there is an associated
filtration
\begin{equation}\label{filt}
\fg=\fg^{-k}\supset\fg^{-k+1}\supset\dots\supset
\fg^0=\fp\supset\fg^1=\fp_+\supset\fg^2\dots\supset\fg^k=\fg_k
\end{equation}
defined by $\fg^j:=\oplus_{\ell\geq j}\fg_\ell$ which is invariant under the adjoint
action of $P$. The subgroup $G_0\subset P$ of all elements whose adjoint action does
preserve the grading has Lie algebra $\fg_0$. This is the reductive Levi factor of
$P$. We also write $P_+=\operatorname{exp}\fp_+$ for the nilradical of $P$, and
similarly $G_- = \operatorname{exp}\fg_-$.  In particular, for each $g\in P$, there
are unique elements $g_0\in G_0$, $\Up\in \fp_+$, and $\Up_i\in \fg_i$,
$i=1,\dots,k$, such that
\begin{equation}\label{aff_structure}
g=g_0\cdot\operatorname{exp}\Up = g_0\cdot\operatorname{exp}\Up_1\cdot\dots\cdot\operatorname{exp}\Up_k
.\end{equation}
This decomposition reflects the fixed splitting of the filtration of
$\fp_+$ by $P$-submodules, i.e. our fixed 
isomorphism $\operatorname{gr}\fp_+\to \fp_+$. 

The Cartan connection $\om$ provides
the {\em constant vector fields} $\om^{-1}(X)\in \frak X(\cg)$ 
defined for all $u\in\cg$ and $X\in\fg$
by
$$
\om(\om^{-1}(X)(u)) = X 
.$$
These generalize the  left invariant vector fields on the homogeneous model $G\to
G/P$ and enjoy the same equivariancy property, i.e.
$$
Tr^g\cdot \om^{-1}(X)(u) = \om^{-1}(\operatorname{Ad}_{g^{-1}} \cdot X)(u\cdot g),
$$
where $r^g$ is the principal right action by the element $g\in P$. In particular, due
to our fixed splitting of $\fg$, there are the {\em horizontal vector fields}
$\om^{-1}(X)$ with $X\in \fg_-$.

On the homogeneous model, the Maurer-Cartan equation reads as $d\om
+\frac12[\om,\om] = 0$, while on a general geometry the same expression provides the two-form
\begin{equation}\label{cartan_curvature}
K = d\om + \frac12[\om,\om] 
\end{equation}
called the {\em curvature}. 
The equivariance properties of the Cartan connection imply that $K$ is
always a horizontal two-form and so the curvature is completely determined
by the {\em curvature function} $\tilde\ka\in C^\infty(\cg, \La^2\fg_-^*
\otimes\fg)$, 
$$
\tilde\ka(u)(X,Y) = K(\om^{-1}(X)(u),\om^{-1}(Y)(u))=
[X,Y]-\om(u)([\om^{-1}(X),\om^{-1}(Y)])
.$$

Of course, the values of both the connection and the curvature function split
according to the corresponding splitting of $\fg$ into
\begin{equation}\label{connection_splitting}
\om=\om_- + \om_0 + \om_+,\quad
\tilde\ka = \tilde\ka_- + \tilde\ka_0 + \tilde\ka_+,
\end{equation} 
but the individual components are not $P$-equivariant.  For instance,
$\ka_-$ is a well defined object only if its values are considered as elements of the
quotient $\fg/\fp$. The latter component of the curvature is called the {\em torsion}
of the Cartan connection.

Clearly, the curvature is the obstruction to the integrability of the
horizontal distribution $\om^{-1}(\fg_-)$ in $T\cg$ and the Cartan
connection is locally isomorphic to its homogeneous model if and only if the
curvature vanishes, see e.g., \cite[Section 1.5]{parabook}. 

\subsection{Natural bundles}\label{nat_bundles}
Every $P$-representation on a vector space $\bv$ provides the homogeneous vector
bundle $G\x_P \bv\to G/P$ and, more generally, the associated vector bundles
$$
\cv M = \cg \x_P\bv \to M
$$
with standard fiber $\bv$ over all manifolds with a
parabolic geometry of the type $G/P$. Shortly, we shall talk about
$P$-modules $\bv$ and the induced {\em natural bundles}. 
In the sequel, we shall restrict ourselves to $P$-modules with a
diagonalizable action of the center of $G_0$. 

For instance, the Cartan connection $\om$ on $p:\mathcal G\to M$ identifies the
tangent bundle $TM=\mathcal G\x_P \fg/\fp=\mathcal G \x_P \fg_-$, $(u,X)\mapsto
Tp(\om^{-1}(X)(u))$. Similarly, $T^*M=\mathcal G\x_P \fp_+$ and the duality is
expressed by the Cartan-Killing form on $\fg$. Observe that for $i<0$, we have the
$P$-invariant subspace $\fg^{i}/\fp\subset\fg/\fp$ induced by the filtration of $\fg$
from \eqref{filt}. This determines a smooth subbundle $T^iM\subset TM$, so we get a
filtration of $TM$. Likewise, the subspaces $\fg^i\subset\fp_+$ for $i>0$ give rise
to a filtration of the cotangent bundle $T^*M$ by smooth subbundles.

A special class of natural bundles is induced by the $G$-modules viewed
as $P$-modules. They are called the {\em tractor bundles}, see \cite{BEG},
\cite{parabook} 
for detailed description in many specific geometries 
and historical links. A very special role is reserved for the
{\em adjoint tractor bundles} $\mathcal A$ coming from the adjoint representation of $G$
on $\fg$. Of course, the Lie bracket itself is
$\operatorname{Ad}$-equivariant and thus there is the algebraic bracket $\{\
,\ \}$ on the adjoint tractors.   

It is well known that sections $\si$ of the induced bundles $\cv M$ are in bijective
correspondence with smooth functions $\tilde \si:\cg\to\bv$ which are $P$-equivariant
in the sense that
$$
s(u\cdot g) = g^{-1}\bullet s(u).
$$ Here $u\cdot g$ denotes the principal right action of $g$ on $u$, while the bullet
denotes the $P$-action on $\bv$.  We have already seen the curvature function
$\tilde\ka:\cg\to \La^2(\fg/\fp)^*\otimes \fg=\La^2\fp_+\otimes\fg$, representing a
section of the adjoint tractor valued form $\ka\in\Om^2(M;\ca)$,
cf. \eqref{cartan_curvature}, which is the curvature of $\om$ viewed as a two-form on
$M$. Of course, the adjoint tractor bundles inherit all the $P$-invariant objects
from $\fg$, including the metric defined by the Cartan-Killing form. Further, the
1-forms on $M$ live in the invariant subbundle $\ca^1$ corresponding to $\fp_+$,
while the vector fields on $M$ can be viewed (with the help of the Cartan connection
on $\cg$ as noticed above) as sections of the quotient $\ca/\ca^0$,
  where $\ca^0=\cg\x_P\mathfrak p$. In particular,
the torsion $\ka_-$ of the Cartan connection is a vector valued two-form on $M$,
exactly as we are used to see the torsions of affine connections on manifolds.

\subsection{The bundle of Weyl structures}\label{weyl-rho}
We briefly remind the impacts of reductions of the parabolic structure groups to
their reductive parts. The reader can find more details in \cite{parabook} and
\cite{CapMettler}.

Given the Cartan geometry $(\Cal G\to M,\om)$, we can form the quotient
  $\Cal G_0:=\Cal G/P_+$ which is a principal bundle over $M$ with structure group
  $P/P_+\cong G_0$. Each reduction of $\Cal G$ to the structure group $G_0\subset P$
  can then be seen as a $G_0$ equivariant smooth section $\bar s:\Cal G_0\to\Cal G$
  of the quotient projection. It is well known that such reductions are in bijective
correspondence with the smooth sections $s$ of the bundle $\pi:\mathcal G/G_0\to M$,
which in our case is equal to $\pi:A=\mathcal G\x_P P/G_0\to M$.  We call $\pi:A\to
M$ the {\em bundle of Weyl structures}, see \cite{Herzlich} and \cite{CapMettler}. 

In particular, we see that Weyl structures always exist globally, and for two such
sections $\hat s$ and $s$, there always exist unique equivariant functions
$\Up:\cg_0\to \fp_+$ and $\Up_i:\cg_0\to \fg_i$ for $i=1,\dots,k$, such
that for all $v\in \cg_0$,
$$
\bar{\hat s}(v) = \bar s(v)\cdot \operatorname{exp}\Up(v) = \bar s(v)\cdot
\operatorname{exp}\Up_1(v)\cdots \operatorname{exp}\Up_k(v)
.$$

Now we come to several crucial observations. 

{\em First}, $\mathcal G\to A$ is a principal fibre bundle with the structure group
$G_0$ and the tangent bundle is the associated bundle $TA = \mathcal G\x_{G_0}
(\fg_-\oplus\fp_+)$ via the adjoint action of $G_0$ on $\fg$. In particular, $TA$
naturally splits into two components $TA=L^- \oplus L^+$ corresponding to the
$G_0$-invariant components $\fg_-$ and $\fp_+$ in $\fg$. Moreover, $L^+$ is the
vertical bundle of $\pi:A\to M$, while $L^-=\pi^*TM$.

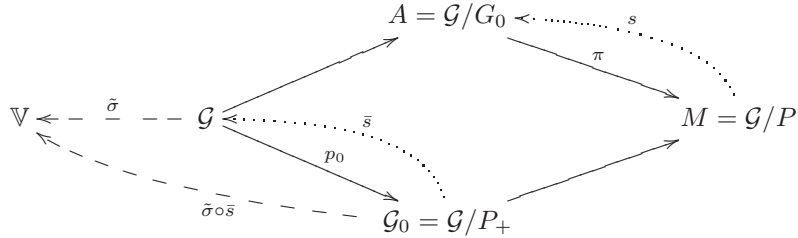
\begin{equation*}
\xymatrix@R=8mm@C=20mm{
&& {A=\mathcal G/G_0} \ar[dr]^{\pi} 
\\
{\bv}
&{\Cal G} \ar[ur] \ar[dr]^{p_0} \ar@(l,r)@{-->}[l]_{\tilde\si} 
&& {M = \mathcal G/P} \ar@(u,r)@{.>}[lu]_{s}
\\
&& {\mathcal G_0 = \mathcal G/P_+} \ar[ur] 
\ar@(u,r)@{.>}[ul]_{\bar s} \ar@(l,dr)@{-->}[llu]^{\tilde\si\o\bar s}
}
\end{equation*}

{\em Second}, the Cartan connection $\om=\om_-+\om_0+\om_+$ can be viewed as an
affine connection $D$ on $A$, with soldering form $\om_-+\om_+$ and principal
connection form $\om_0$. The natural splitting $TA=L^-\oplus L^+$ leads to two
partial affine connections $D^-$ and $D^+$ by restriction. We also observe that,
up to universal algebraic terms, the curvature of $\omega$ on $A$
encodes the torsion and curvature of $D$ by means of the $G_0$-equivariant functions
$\tilde\ka_-+\tilde\ka_+$ and $\tilde\ka_0$ from \eqref{connection_splitting},
respectively. Moreover, since $\ka$ is horizontal with respect to $\mathcal G\to M$,
the torsion and curvature of $D$ are given by universal algebraic terms
  if evaluated on one argument from $L^+$. Evaluated on two arguments from
  $L^-$, they provide, after algebraic corrections, the components $T\in
\Om^2(A,L_-)$, $Y\in\Om^2(A,L^+)$, and $W\in\Om^2(A,\Cal G\x_P\fg_0)$. We call these
components the {\em universal torsion} $T$, {\em universal Cotton York tensor} $Y$,
and {\em universal Weyl curvature} $W$, respectively, see
  \cite{CapMettler} for details.

For each Weyl structure $s$, we can also consider the pullback of the Cartan
connection $\bar s^*\om\in\Om^1(\mathcal G_0,\fg)$ which now splits into the 
three components naturally. 
Let us write
\begin{align*}
\th^s & =\bar s^*\om_- = \th^s_{-k}+\dots+\th^s_{-1}:T\cg_0\to
\fg_{-k}\oplus\dots\oplus\fg_{-1}
\\
\ga^s &= \bar s^*\om_0: T\cg_0\to \fg_0
\\
\Rho^s &= \bar s^*\om_+ = \Rho^s_1+\dots +\Rho^s_k:T\cg_0\to \fp_+ =
\fg_1\oplus\dots\oplus\fg_k 
.\end{align*}
The filtration on $TM$ introduced in Section \ref{nat_bundles} gives rise to the associated
graded bundle $\gr(TM)=\oplus_{i=-k}^{-1}T^iM/T^{i+1}M$. Now the form $\th^s$
provides an identification of $\cg_0$ with a frame bundle for $\gr(TM)$ with
structure group $G_0$ and should be understood as an analog of the {\em soldering
  form}. In particular, this defines the isomorphism $TM\simeq \gr(TM)$, i.e. a
splitting of the filtration
$$
TM=T^{-k}M\supset \dots\supset T^{-1}M
$$
determined by the choice of $s$. 

The next component $\ga^s$ is the
connection form of a principal connection on $\cg_0$. Thus, the first two
components provide a reductive Cartan connection $(\th^s+\ga^s)$ of type 
$(G_-\cdot G_0)/G_0$ on $M$. We call it the {\em Weyl connection}
corresponding to
the choice of the Weyl structure $s$. We shall write $\nabla^s$ for the
covariant derivatives with respect to this connection on all vector bundles
associated to $\cg_0$.

The last component $\Rho^s$ measures the difference
between the Cartan connection $(\th^s+\ga^s)$ on $\cg_0$ and the original Cartan
connection $\om$ on $\cg$, along the image of $\bar s$. 
We view it as a one-form on $M$ valued in
$\operatorname{gr}T^*M$ and call it the {\em Rho-tensor}.

{\em Third}, let us take a representation $\bv$ of $P$ and consider the
  induced bundle $\Cal G\x_P\bv=\cv M\to M$. Sections of this bundle correspond to
  $P$-equivariant maps $\tilde \si:\Cal P\to\bv$, so they obviously also represent
  sections of the natural bundle $\mathcal V A = \mathcal G\x_{G_0} \bv\to A$. By
the very construction, the natural bundles $\cv A\to A$ are naturally identified with
$\pi^*\cv M$.

Next, the action of the grading element $E\in\fg_0$ decomposes $\bv$ into
$G_0$-invariant components
$$
\bv = \bv_0\oplus\dots\oplus\bv_\ell
$$
with the property that $X\bullet\bv_j\subset \bv_{i+j}$, for all $X\in\fg_i$ and 
$i=1,\dots,k$. For instance, the Lie algebra $\fg_-$, viewed as
$\operatorname{gr}(\fg/\fp)$, or the entire $\fg$ are special cases
illustrating the difference between the $P$-invariant filtrations and
$G_0$-invariant induced gradings. 

Thus while the natural bundles $\cv M\to M$ come equipped with the natural
filtrations only, the bundles $\cv A \to A$ come with natural gradings.  Once we
choose a Weyl structure $s$, the bundles $\cv M$ get graded, too (via the structure
group reduction). This grading coincides with the pullback $s^*(\cv A)\to M$,
\cite[Theorem 2.4]{CapMettler}.

Of course, the canonical affine connection $D$ induces linear connections on all
natural bundles $\cv A$.  Now, the required $P$-equivariance of the maps $\tilde \si$
representing sections of $\cv M$ directly implies that the sections of $\cv M$ are
exactly those sections $\tau$ of $\cv A$, whose covariant derivatives satisfy $D_\xi
\si=-\xi\bullet\tau$ for all sections $\xi\in \Ga(L^+)$, \cite[Theorem
  2.4]{CapMettler}.

For a section $\si\in\Ga(\cv M)$, the covariant derivative $D\si$ with
  respect to the canonical connection $D$ defines a map $D \si:\cg \to (\fg_-\oplus
  \fp_+)^*\otimes\bv$. Restricting to entries from $\fg_-$, one obtains an operation
  $\nabla^\om$ which is called the {\em invariant differential} and which played a
  crucial role in \cite{CSS1-3} and \cite{CSS4}. However, this was viewed as an
  operation mapping sections of $\cv M$ to equivariant functions on the Cartan
  bundle, the relation to the bundle of Weyl structures was not known at that time. 

\subsection{The curvatures}\label{curvature_comparison}
As mentioned in Section \ref{cartan-connection}, the fundamental invariant of a
Cartan geometry $(\Cal G,\om)$ is the Cartan curvature $\ka$. Via the curvature
function discussed there, this admits direct interpretations both on $M$ and on
$A$. $P$-equivariancy of the curvature functions implies that $\ka$ can be viewed as
a two-form on $M$ with values in the adjoint tractor bundle $\Cal A=\Cal G\x_P\frak
g$. A choice of Weyl structure determines an isomorphism between the bundle
$\La^2T^*M\otimes\Cal A$ with its associated graded. The latter bundle decomposes
into components according to the decomposition of $\La^2(\frak g_{-})^*\otimes\frak
g$ into $G_0$-irreducible representations.

On the level of $A$, the decomposition into $G_0$-components is available
canonically. In particular, the $G_0$-invariant decomposition $\fg=\fg_-\oplus
\fg_+\oplus\fg_0$ decomposes $\Cal G\x_{G_0}\fg$ with the first two summands
corresponding to $L^-\oplus L^+=TA$. Via the restriction of the adjoint action,
$\fg_0$ canonically injects into endomorphisms of $\fg_-\oplus\fg_+$, so we can view
$\Cal G\x_{G_0}\fg_0$ naturally as a subbundle $\End_0(TA)\subset T^*A\otimes TA$. So
on this level, we can interpret $\ka$ as being given by a section of the bundle
$\La^2(L^-)^*\otimes TA$ which we write as $T+Y$ according to $TA=L^-\oplus L^+$ and
a section $W\in\Ga(\La^2(L^-)^*\otimes\End_0(TA))$. These look already like a torsion
and a curvature, and indeed by \cite[Theorem 2.12]{CapMettler} these equivalently
encode the torsion and curvature of $D$. More precisely, one has to add universal
bundle maps induced by ($G_0$-equivariant) components of the Lie bracket on $\fg$ to
$T+Y$ respectively to $W$ in order to obtain the torsion and the curvature of $D$.

Next, we want to relate the Cartan curvature to the curvature and torsion of the Weyl
connections associated to a Weyl structure. Consider a Weyl structure $s$, and the
corresponding affine connection $\th^s+\ga^s$ on $\cg_0$. The full Cartan curvature
(i.e., torsion and curvature) of this connection is also given by the structure
equations
\begin{align*}
T^s(\xi,\eta) &= d\th^s(\xi,\eta) + [\ga^s(\xi),\th^s(\eta)]) +
[\th^s(\xi),\ga^s(\eta)]
+ [\th^s(\xi),\th^s(\eta)]
\\
R^s(\xi,\eta) &= d\ga^s(\xi,\eta) + [\ga^s(\xi),\ga^s(\eta)],
\end{align*}
where $\xi,\eta$ are vector fields on $\cg_0$. Clearly the torsion $T^s$ and the
curvature $R^s$ are horizontal two-forms which descend to well defined forms on the
underlying manifold. Let us point out, that $T^s$ is not the usual torsion of a
linear connection but there is a correction involving natural bundle maps induced by
the ($G_0$-equivariant) Lie bracket on $\fg_-$.  Since the forms $\th^s$ and $\ga^s$
are pullbacks of $\om_-$ and $\om_0$, the above structure equations are easily
compared with those of $\om$.

The latter curvature forms are clearly related to the pullback $\bar s^*\ka$ of the
curvature of the Cartan connection $\om$ on $\cg$ (i.e.  also the torsion and
curvature of $D$).  The missing components of curvature are related to the
$\Rho$-tensor and they can be understood easily, too:
\begin{equation}\label{WeylCurvature}
\bar s^*\ka = T^s+R^s+Y^s + \partial \Rho^s
\end{equation}
where we define the {\em Cotton-York} tensor of the Weyl structure on $M$ as
$$
Y^s(\xi,\eta) = d^{\nabla^s}\Rho^s(\xi,\eta) + \Rho^s(\{\xi,\eta\}) +
\{\Rho^s(\xi),\Rho^s(\eta)\}
,$$
with $d^{\nabla^s}$ denoting the covariant exterior differential, $\{\ ,
\ \}$ is the natural bracket on the adjoint tractors, 
and 
$$\partial \Rho^s(\xi,\eta) = \{\xi,\Rho^s(\eta)\} - \{\eta,\Rho^s(\xi)\} -
\Rho^s(\{\xi,\eta\})
$$ is the ($G_0$-equivariant) Lie algebra cohomology differential. One then defines
the \textit{Weyl curvature} of the Weyl connection determined by $s$ as
$W^s:=R^s+(\partial\Rho^s)_{\mathfrak g_0}$.

As suggested by the notation, the quantities $T^s$, $W^s$ and $Y^s$ are directly
related to the quantities on $A$ discussed above. Indeed by \cite[Proposition
  2.10]{CapMettler} and \cite[Section 5.2.9]{parabook} one has $s^*T = T^s+(\partial
\Rho^s)_{\mathfrak g_-}$, $s^*Y = Y^s+(\partial \Rho)_{\mathfrak p_+}$, and
$s^*W=W^s$. There is also a nice way to interpret the Rho-tensor via $A$. Viewing the
projection onto the second factor in $TA=L^-\oplus L^+$ as an element of
$\Rho\in\Om^1(A,L^+)$, one can form $s^*\Rho\in\Om^1(TM,\gr(TM))$ and by
\cite[Proposition 2.8]{CapMettler} this coincides with $\Rho^s$.

It will be very important for our results that the decomposition of the
  curvature and torsion of a Weyl connection into the pullback of the Cartan
  curvature and a part obtained from the Rho-tensor can be made explicit without
  reference to the Cartan bundle. This only works for \textit{normal parabolic
    geometries} which anyway is the class usually studied. They are characterized by
  a normalization condition on their curvature $\ka$ which is usually phrased as
  $\partial^*(\ka)=0$. Here $\partial^*$ is induced by a $P$-equivariant linear map
  $\La^2(\fg/\frak p)^*\otimes\fg\to (\fg/\frak p)^*\otimes\fg$, the so-called
  \textit{Kostant codifferential}. This normalization condition then implies that the
  parabolic geometry is determined by some underlying geometric structure known as an
  \textit{infinitesimal flag structure}, for example a conformal structure.

  The detailed form of this is not important for our purposes. What it tells us,
  however, is that the right hand side of \eqref{WeylCurvature} lies in the kernel of
  $\partial^*$. This allows us to express, $\partial^*\partial(\Rho^s)$ in terms of
  the curvature and torsion of $\nabla^s$ and of $Y^s$ as defined above. Now on the
  one hand it turns out that $\partial^*\partial$ is invertible. On the other hand,
  the definition of $Y^s$ is chosen in such a way our equation can be solved
  iteratively homogeneity by homogeneity (with respect to the grading element), see
  Section 5.1 of \cite{parabook} for details. Hence $\Rho^s$ can be computed from
  curvature quantities, which in turn allows us compute the pullback of the Cartan
  curvature. In particular, in the case of conformal geometry this leads to the
  description of the Rho tensor as a trace-modification of the Ricci curvature as
  introduced by J.A.\ Schouten.

\subsection{Normal Weyl structures}\label{normal_weyl}
As we have seen, the Rho-tensor measures the deviation of the Cartan connection
$\th^s+\ga^s$ on $\cg_0$ determined by a Weyl structure $s$, from the given Cartan
connection $\om$. Thus minimizing the values of $\Rho^s$ and its derivatives in a
point looks like a good idea. We can easily follow the way how the normal coordinates
for affine connections are built. This leads to an analog for all parabolic
geometries of the so-called Graham normal coordinates (see \cite{FefGr}) from
conformal geometry and their analog in CR geometry that was recently introduced by
K. Hirachi in \cite{Hi}.

For each fixed frame $u\in \cg$ and $X\in\fg_-$, there is the flow line
$\operatorname{Fl}^{\om^{-1}(X)}_t(u)$ of the horizontal vector field $\om^{-1}(X)$,
and for some neighborhood of the origin in $\fg_-$ these flow lines exist at least up
to the time $t=1$. In this way we obtain a horizontal embedding $X\mapsto
\ph_u(X)=\operatorname{Fl}^{\om^{-1}(X)}_1(u)$ and also a local section
$p(\ph_u(X))\mapsto \ph_u(X)$ of the Cartan bundle $p:\cg\to M$ through $u$.
Consequently there is a unique local Weyl structure $\bar s_u$ through the frame $u$
defined by $\bar s_u(p_0(\ph_u(X))) = \ph_u(X)$. (Here $p_0:\cg\to \cg_0$ is the
natural projection and the $G_0$-orbit of this section defines the reduction.)  We
call them the {\em normal Weyl structures}. The images $c^{u,X}$ of the defining flow
lines in $M$, i.e. $c^{u,X}(t) = p(\operatorname{Fl}^{\om^{-1}(X)}_t)$, are called
the {\em generalized geodesics} of the Cartan connection $\om$. Thus, each choice of
a frame $u\in\cg$ over the point $p(u)=x\in M$ determines a uniquely defined
``geodetical parametrization'' from a neighborhood of $0$ in $T_xM$ to a neighborhood
of $x$ in $M$.

This construction can be also nicely seen on $A$. We are restricting the standard
normal coordinates of the affine connection $D$ to geodesics emanating in
$L^-$-directions, which directly creates a (local) section $s:M\to A$, similarly as
above. These sections evidently have the property, that the universal Rho tensor
vanishes along those geodesics. Indeed, this leads to a complete characterization:
The Weyl structures $s$ that are normal at $x$ are characterized by the
  fact that $\Rho^s(c(t))(c'(t))=0$ for all generalized geodesics through $x=p(u)\in
  M$.  In particular, this implies that for every integer $\ell\ge0$ of covariant
derivatives, and any tangent vectors $\xi_0,\dots,\xi_\ell$, the full symmetrization
of the expression
\begin{equation}
(\nabla^s_{\xi_\ell}\cdots\nabla^s_{\xi_1}\Rho(x))(\xi_0)\in \Om^1(M)
\end{equation}
over the $\xi$'s vanishes, see \cite{CS} for details. 

\section{The invariant calculus}

Whatever concept of invariance we adopt, all objects built naturally in the category
of Cartan connections of given type should be invariant. Thus, we can consider any
(possibly non-linear) differential operator $\Ph:\Ga(\cv A) \to \Ga(\mathcal W A)$
between two natural vector bundles over $A$ coming from $P$-modules $\bv$ and
$\mathbb W$, expressed in covariant derivatives with respect to $D$ and derivatives
of the torsion and curvature of $D$. As noted in Section \ref{weyl-rho}, we can view
$\Ga(\cv M)$ naturally as a subspace of $\Ga(\cv A)$ (corresponding to
$P$-equivariant functions among $G_0$-equivariant functions). Then it may happen that
the restriction of $\Ph$ to this subspace has values in $\Ga(\mathcal W M)\subset
\Ga(\mathcal W A)$. If this is the case, the restriction has to define an invariant
operator.

However, choosing a Weyl structure $s$, there is a possibility to
``descend'' any operator $\Ph:\Ga(\cv A) \to \Ga(\mathcal W A)$ as above to
an operator $\Phi^s:\Ga(\cv M) \to \Ga(\mathcal W M)$.  Namely, for a
section $\si\in\Ga(\cv M)$ we can view $\si$ as an element of $\Ga(\cv A)$
and form $\Phi(\si)\in\Ga(\mathcal W A)$.  Then we can take the pullback
$s^*\Phi(\si)$ and observe that there is a unique section
$\Phi^s(\si)\in\Ga( \mathcal W M)\subset\Ga(\mathcal W A)$ such that
$s^*\Phi^s(\si)=s^*\Phi(\si)$.  Equivalently, this can be phrased as
$\Phi(\si)|_{s(M)}=\Phi^s(\si)|_{s(M)}$.  In the language of equivariant
functions, this simply means that we take the $G_0$-equivariant function
$\Phi(\si):\cg\to\mathbb W$, restrict it to $\tilde s(\cg_0)\subset\cg$ and
then define $\Phi^s(\si):\cg\to\mathbb W$ to be the $P$-equivariant
extension of this $G_0$-equivariant function.

Observe that this construction immediately implies a property that will be
of crucial importance in what follows, namely that for any $\si\in\Ga(\cv M)$
and any point $x\in M$, the value $\Phi^s(\si)(x)$ depends only on $s(x)$,
so the dependence of $\Phi^s$ on $s$ is of tensorial character.

An obvious first step in trying to understand this construction is to
analyze the operator $\Phi:\Ga(\cv A)\to \Ga((L^-)^*\otimes\cv A)$ defined
by $\Phi(\si):=D^-\si$ on some natural bundle.  Choosing a Weyl structure
$s$, this descends to an operator $(D^-)^s:\Ga(\cv M)\to \Ga(T^*M\otimes\cv
M)$ which has the chance to define a covariant derivative on $\cv M$.  This
case is sorted out by a result which basically is proved as a part of
Theorem 2.6 in \cite{CapMettler}:

\begin{prop}\label{nabla_comparison}
  Let $\bv$ be any representation of $P$ and consider the corresponding natural
  bundles $\cv M\to M$ and $\cv A\to A$, and the operator $\Phi:\Ga(\cv A)\to
  \Ga((L^-)^*\otimes\cv A)$ defined by $\Phi(\si):=D^-\si$. Then for any Weyl
  structure $s:M\to A$, the operator $\Phi^s:\Ga(\cv M)\to\Ga(T^*M\otimes\cv M)$
  coincides with the Rho-corrected derivative $\nabla^{\smallRho^s}$ corresponding to
  $s$ as introduced in Section 5.1.9 of \cite{parabook}. 
\end{prop}
\begin{proof}
Consider $\bar s(\cg_0)\subset\cg$ and restrict the component $\omega_\fp$ to this
subset. Then this can be uniquely extended to a principal connection form $\ga^s$ on
$\cg$ and by definition, the Rho-corrected derivative is the covariant derivative
induced by this principal connection. Viewing $\si\in\Ga(\cv M)$ as a $P$-equivariant
function $\cg\to\bv$, the function $\cg\to L(\fg_-,\bv)$ corresponding to
$\nabla^{\smallRho^s}\si$ thus is given by taking the derivative of $\si$ with
respect to horizontal vector fields along $\bar s(\cg_0)$ and then equivariantly
extending to $\cg$. But this exactly means that, along $\bar s(\cg_0)$, this
coincides with the $G_0$-equivariant function associated to $D^-\si$, see Section
\ref{weyl-rho}, which proves the claim. 
\end{proof}

In view of our observation above, this gives an alternative argument for the fact
that the transformation rule for Rho-corrected derivatives is tensorial in the
one-form $\Upsilon$ describing the change between two Weyl-structures, see part (2)
of \cite[Proposition 5.1.9]{parabook}.

The name Rho-corrected derivative is motivated by the relation of this operation
(determined by some Weyl structure $s$) to the Weyl connection $\nabla^s$.  This can
be written as $\nabla^{\smallRho^s}_\xi\si=\nabla^s_\xi\si+\Rho^s(\xi)\bullet\si$ for
any vector field $\xi\in\frak X(M)$, since the difference of the horizontal lifts of
$\xi\in T_xM$ for $\ga^s$ and $\ga$ on the image $\bar s(\mathcal G_0)$ is the
fundamental vector field corresponding to $-\om_{\mathfrak p_+}(\xi)$, and $\si$ is
$P$-equivariant. Here the bullet comes from the bundle map $T^*M\x \cv M\to\cv M$
induced by the infinitesimal representation $\frak p_+\x\bv\to\bv$.  A bit of care is
needed in interpreting this formula, however, since both the Weyl connection
$\nabla^s$ and the $\Rho^s(\xi)$ are actually carried over from the associated graded
bundles $\gr(\cv M)$ and $\gr(T^*M)$ via the isomorphism induced by $s$.  This also
explains how the transformation law under a change of Weyl structure can be tensorial
in the one-form $\Upsilon$ in spite of the explicit occurrence of a Rho-tensor (whose
transformation is not tensorial, see Section \ref{transf} below).  Hence under the
isomorphism $\Ga(\Cal VM)\cong \Ga(\gr(\Cal VM))$ induced by a Weyl structure $s$,
the operator $(D^-)^s$ corresponds to a universal expression in terms of $\nabla^s$
and the Rho-tensor $\Rho^s$ of $s$.

We next want to obtain similar descriptions for more general operators built up from
$D^-$ and the torsion and curvature of $D$. The first step is to understand the
iterated derivatives $(D^-)^k$ with respect to $D^-$. This is tricky, because the
argument in the proof of Proposition \ref{nabla_comparison} was based on the fact
that we start with $\si\in\Ga(\cv M)$, i.e.\ with a $P$-equivariant
function. However, the function corresponding to $D^-\si$ is not $P$-equivariant any
more, so we cannot simply iterate. A neat way around this is via an invariant notion
of jets that was crucial for the developments in \cite{CSS1-3} and in \cite{CSS4}.

\subsection{First order jets}\label{1-jets}
In geometry, a differential operator $\Ga(\cv M)\to\Ga(\Cal WM)$ of
  order at most $k$ can be equivalently viewed as a vector bundle homomorphism $J^k\cv
  M\to \mathcal W M$, where $J^k\cv M$ is the $k$th jet prolongation of $\Cal VM$.
Thus, in order to understand the operators on $\cv A$ and $\cv M$, let us look at 
jets.
 
Let us first consider Klein geometries, i.e., the flat model Cartan geometries $G\to
G/P$ with the Maurer-Cartan form $\om$. There, natural vector bundles are exactly the
homogenous bundles $\cv = G\x_P\bv$ determined by $P$-modules $\bv$. Any jet
prolongation of a homogeneous bundle is homogeneous, too, so there is a $P$-module
$J^k\bv$ such that $J^k\cv=G\x_P J^k\bv$. The $P$-module $J^k\bv$ can be obtained as
the fiber of $J^k\cv$ over the origin $o\in G/P$.  Thus, the invariant operators are
in bijective correspondence with the intertwining maps $J^k\bv\to\mathbb W$ between
the $P$-modules. In a dual picture, this provides an identification of linear
invariant operators with the singular vectors in generalized Verma modules for all
parabolic models, see e.g. \cite{ES}.

Let us now consider a parabolic geometry $\cg\to M$ with Cartan connection
$\om$. In \ref{weyl-rho}, we identified the covariant derivative of the canonical
affine connection $D$ restricted to $L^-$ with the so called invariant differential
$\nabla^\om$ acting on $\bv$-valued functions on $\cg$ via the constant vector fields
$\om^{-1}(X)$, $X\in\fg_-$.  This is just the usual definition of the covariant
derivative of affine connections in its frame form.

Although $\nabla^\om$ does not map $P$-equivariant functions $\si$ into
$P$-equivariant results, its extension
$$
D^\om : C^\infty(\cg,\bv) \to C^\infty(\cg, \fg^*\otimes\bv),\quad 
D^\om \si(u)(X) = \om^{-1}(X)(u)\cdot\si \in \bv
$$
does. This operator is called the {\em fundamental derivative}. It is very well
understood in terms of the adjoint tractor bundles and leads to the so called tractor
calculus, see \cite[Sections 1.5.7-8]{parabook}.

The operator $D^\om:\Ga(\cv M)\to \Ga(\mathcal W M)$ can be understood as the
universal differential operator and the standard fiber $\Bbb W$ of the target is $\fg^*\otimes
\Bbb V$. We may also consider the analogue of the first jet prolongation by mapping
$\si$ to $(\si,D^\om\si)$. The natural action of $g\in P$ on this target module $\Bbb
V\oplus(\fg^*\otimes \Bbb V)$ is
$$
g\cdot(v,\ph) = (g\cdot v, B\mapsto g\cdot\ph(\operatorname{Ad}_{g^{-1}}B)) .
$$
The definition of $D^\om$ together with equivariancy of the function $\si$ that is
differentiated readily implies that for $X\in\frak p\subset\frak g$, we get
$D^\om\si(u)(X)=-X\bullet\si(u)$, where in the right hand side, we use the
infinitesimal action of $\frak p$ on $\bv$. Correspondingly, there is a $P$-submodule
$J^1\bv\subset\bv\oplus \fg^*\otimes\bv$, and $\si\mapsto (\si,D^\om\si)$ induces a
functorial isomorphism $J^1\cv M\to \Cal G\x_P J^1\bv$. This construction works for
all Cartan geometries of type $G/P$.  See \cite[Section 1.5.10]{parabook} for
details.

Now as a $G_0$-module, $\fg=\fg_-\oplus\fp$ which leads to an isomorphism
$J^1\bv\cong_{G_0} \bv\oplus(\fg_-)^*\otimes\bv$ via restriction in the second
component. Consequently, we can naturally identify $\Cal G\x_{G_0}J^1\bv\to A$ with
$\cv A\oplus (L^-)^*\otimes\cv A$. Now for a section $\si\in\Ga(\cv M)$ the value of
the one-jet operator $j^1\si\in\Ga(J^1\cv M)$ defines a section of $\Cal
G\x_{G_0}J^1\bv\to A$, which under this identification corresponds to $(\si,D^-\si)$. 

\subsection{Higher order jets}
The iteration of the first jet prolongation leads to the non-holonomic jet
prolongations and most of the redundancies can be removed by considering
semi-holonomic jets.  It turns out that the target space is a natural bundle with the
fiber called {\em semiholonomic jet module} $\bar J^k\bv$ induced by $\bv$. Hence
there is a $k$th order {\em invariant semiholonomic jet} operator, $k=1,2,\dots$ with
values in $\Cal G\x_P\bar J^k\bv$, see \cite{CSS1-3} for details. Now similarly as in
order 1, we get an isomorphism of $G_0$-modules
\begin{equation}\label{jet-iso}
  \bar J^k\bv\cong_{G_0}\bv\oplus ( \fg_-)^*\otimes\bv\oplus\dots\oplus\otimes^k
  ( \fg_-)^*\otimes\bv
\end{equation}
which is induced by restriction of multilinear maps. In particular, since any bundle
map induced by a $G_0$-homomorphism is parallel for $D$ and hence for $D^-$, for a
section $\si\in\Ga(\cv M)$, the semi-holonomic $k$-jet operator $\bar j^k\si$
corresponds to $(\si,D^-\si,\dots,(D^-)^k\si)$ under this identification. We will
implement this description by viewing $ \bar J^k\bv$ as the $G_0$--module $\bv\oplus
( \fg_-)^*\otimes\bv\oplus\dots\oplus\otimes^k (
  \fg_-)^*\otimes\bv$ endowed with a (very 
complicated) extension of the action to $P$. Still, the infinitesimal action of
$\fp_+$ is given
by a $G_0$--equivariant map $\frak p_+\otimes \bar J^k\bv\to \bar J^k\bv$ in this
picture.

From the construction it follows easily that for $\cv M=\Cal G\x_P\bv$ the invariant
semi-holonomic jet operator defines an injective bundle map $J^k\cv M\hookrightarrow
\Cal G\x_P\bar J^k\bv$. However, the image does not correspond to a $P$-invariant
subspace in $\bar J^k\bv$ and there seems to be no way to realize higher order jet
prolongations as associated bundles to the Cartan bundle on Cartan geometries with
non-vanishing curvature. (The identification of a jet prolongation with
  an induced bundle on the homogeneous model $G/P$ discussed in Section \ref{1-jets}
  extends to geometries with vanishing curvature, which are locally isomorphic to
  $G/P$ but in general not further.)



This makes the construction of invariant differential operators via algebraic
techniques subtle. On the one hand, any $P$-equivariant map $J^k\bv\to \mathbb W$ induces a
$k$th order invariant differential operator on locally flat geometries, but not on
general curved geometries. On the other hand, any $P$-equivariant map $\bar
J^k\bv\to\mathbb W$ gives rise to an invariant differential operator on all Cartan
geometries (via composition with the invariant semi-holonomic jet operator) but, even
in the linear case, not all invariant differential operators are obtained in this
way. One obtains many operators in this way, however, in particular all operators that
arise in BGG-sequences are of this character, see e.g. \cite{ES} and \cite{CSS4}. We
will use the semi-holonomic jet prolongations in a different way here.

\subsection{Expansion of $(D^-)^k$}\label{higher-order}
Now we are ready to prove our first main result, which is a higher order analog of
Proposition \ref{nabla_comparison}, so we want to study the operators
$((D^-)^k)^s$ for
a Weyl structure $s$ and any $k\geq 1$. We will show that, under the isomorphism to
the associated graded, this operator can be written via a ``universal formula'' in
the following sense. Fix an initial representation $\bv$ of $P$ and an order
$k$. Then we need finitely many $G_0$-equivariant linear maps $A_i$ with values in
$\otimes^k\fg_-\otimes\bv$, such that for each Weyl structure $s$ and any section
$\si\in\Ga(\cv M)$, we can write the image of $s^*((D^-)^k(\si))$ in the associated
graded (under the isomorphism determined by $s$) as
$(\nabla^s)^ks^*\si+\sum_i\tilde A_i(T_i)$. Here $\nabla^s$ denotes the Weyl
connection, $\tilde A_i$ is the bundle map induced by $A_i$ and $T_i$ is an iterated
tensor product of factors of the form $(\nabla^s)^\ell s^*\si$ and
$(\nabla^s)^\ell\Rho^s$ with $0\leq\ell<k$ that is independent of $s$. (The form of such
a tensor product also determines the domain of the map $A_i$.)

\begin{thm}\label{k-th-order}
  For any representation $\bv$ and any $k\geq 1$, there is a universal expression for
  the image of $s^*((D^-)^k(\si))$ in the associated graded under the isomorphism
  determined by $s$ as $(\nabla^s)^ks^*\si$ plus a sum of terms
  containing iterated Weyl 
  derivatives of $\si$ and $\Rho^s$ in the sense described above. 
    Moreover, in any tensor product in this expansion, the total number
    of covariant derivatives is $<k$.
\end{thm}
\begin{proof}
  We proceed by induction on $k$. For $k=1$, Proposition \ref{nabla_comparison}
  provides an expansion in the claimed form with $i=1$,
  $A_1:(\fg_-)^*\otimes\fp_+\otimes\bv\to (\fg_-)^*\otimes\bv$, the tensor product of
  the identity with the infinitesimal representation, and $T_1=\Rho^s\otimes\si$.

  Assuming that the result has been proved for all $\ell\leq k$, we consider the
  invariant semi-holonomic $k$-jet operator $\bar j^k$, which via our description of
  $\bar J^k\bv$ is given by $\bar j^k\si=(\si,D^-\si,\dots,(D^-)^k\si)$. Of course,
  the components of the pullback $s^*(\bar j^k\si)$ are just $s^*(D^-)^\ell\si$ for
  $0\leq\ell\leq k$. Hence by induction hypothesis, we can express the image of
  $s^*\bar j^k\si$ in the associated graded under the isomorphism induced by $s$ as
  $(\nabla^s)^ks^*\si$ plus a sum  $\sum_i\tilde A_i(T_i)$ of universal
    terms as described above which involve covariant derivatives of $s^*\si$ and of
    $\Rho^s$. Moreover, the total number of covariant derivatives occurring in each
    $T_i$ is $<k$.

  Since $\bar j^k\si$ lies in the subspace $\Ga(\Cal G\x_P\bar J^k\bv)$, we can apply
  Proposition \ref{nabla_comparison} to conclude that
  $s^*(D^-\bar j^k\si)=\nabla^s(s^*\bar j^k\si)+\Rho^s\bullet s^*\bar j^k\si$. But
  naturality of $D^-$ implies that viewing $(\fg_-)^*\otimes\bar J^k\bv$ as
  $\oplus_{r=1}^{k+1}\otimes^r(\fg_-)^*\otimes\bv$, the components of
  $D^-\bar j^k\si$ are just $(D^-)^r\si$ for $r=1,\dots,k+1$. In particular, we can
  extract an expression for $s^*(D^-)^{k+1}\si$ from $s^*(D^-\bar j^k\si)$ and the
  known expressions for the iterated derivatives of lower order.

  Expanding $s^*\bar j^k\si=(\nabla^s)^ks^*\si+\sum_i\tilde A_i(T_i)$ and applying
  $\nabla^s$, we obtain $(\nabla^s)^{k+1}s^*\si$ from the first summand. On the other
  hand $\Rho^s\bullet (\nabla^s)^ks^*\si$ is of the form that we allow for the
  additional terms in the expansion of $s^*((D^-)^{k+1}\si)$. Thus we can conclude
  the proof by showing that the expressions $\nabla^s(\tilde A_i(T_i))$ and
  $\Rho^s\bullet \tilde A_i(T_i)$ all lead to expressions of the form we allow. For
  the first case, we can use that the bundle map $\tilde A_i$ is induced by a
  $G_0$-equivariant linear map and hence it is parallel for $\nabla^s$. Thus we
  obtain $(\id\otimes\tilde A_i)(\nabla^sT_i)$ and of course $\nabla^sT_i$ can be
  expanded as a tensor product of factors $(\nabla^s)^\ell\si$ and
  $(\nabla^s)^\ell\Rho^s$ with a total of at most $k$ covariant
    derivatives by our assumptions on $T_i$. For the other term, we can write
  $\Rho^s\bullet \tilde A_i(T_i)$ as $\tilde B_i(\Rho^s\otimes T_i)$. Here
  $B_i=(\id\otimes\tau)\o (\id\otimes\id\otimes A_i)$, where $\tau$ denotes the
  infinitesimal representation, so this is again of the required form.
\end{proof}

\subsection{Involving curvature quantities}\label{involve-curv}
As discussed in Section \ref{curvature_comparison}, the torsion and the curvature of
the canonical connection $D$ on $A$ can be equivalently encoded via the curvature
$\ka$ of the canonical Cartan connection $\om$.  Now the Cartan curvature $\ka$ is a
section of $\Cal G\x_P\bv\to M$, where $\bv:=(\La^2(\fg/\fp)^*\otimes\fg)$, so
Theorem \ref{k-th-order} applies to $\ka$.  Hence for any $k\in\Bbb N$, we can
describe the pullback $s^*(D^-)^k\ka$ in terms of $(\nabla^s)^ks^*\ka$ and universal
terms depending on iterated Weyl derivatives of order less than $k$ of $s^*\ka$ and
of the Rho-tensor $\Rho^s$. We have seen, that on $A$, $\ka$ decomposes as $T+W+Y$
according to the values in $\mathfrak g=\mathfrak g_-\oplus\mathfrak
g_0\oplus\mathfrak p_+$. Of course, these and further decompositions coming from
$G_0$-invariant operations are, on the level of $A$, preserved by the operator
$D^-$. They correspond to decompositions of $s^*\ka$ on $M$ (depending on $s$) which
in turn are preserved by the Weyl connection $\nabla^s$.  We also discussed in
Section \ref{curvature_comparison} how to decompose $s^*\ka$ into its components
$s^*T=T^s+(\partial \Rho^s)_{\mathfrak g_-}$, $s^*W = R^s +(\partial
\Rho^s)_{\mathfrak g_0}$, and $s^*Y=Y^s+({\partial \Rho})_{\mathfrak p_+}$, which are
related to the torsion and curvature of $\nabla^s$.

Otherwise put, if we consider a component $K$ of the curvature and torsion of
$D$, and form $(D^-)^kK$, we can recover this from $(D^-)^k\ka$ by $G_0$-equivariant
operations. Applying the same operations to $s^*(D^-)^k\ka$ we conclude that we can
express $s^*(D^-)^kK$ as $(\nabla^s)^ks^*K$ plus a universal expression in the sense
discussed in Section \ref{higher-order} involving iterated derivatives of order $<k$
of components of $s^*\ka$ and of $\Rho^s$.

But this is sufficient to pull back polynomial invariant operators (in a sense that
will be made precise) constructed from $D^-$ along Weyl structures and obtain a
special class of affine invariants of Weyl structures. So we start with arbitrary
representations $\bv$ and $\Bbb W$ of $P$ and we are looking for an operator
$\Ph:\Ga(\Cal VA)\to\Ga(\Cal WA)$ such that we can write $\Ph(\si)=\sum_i A_i(T_i)$,
where $A_i$ is a natural bundle map induced by a $G_0$-equivariant map with values in
$\Bbb W$ and $T_i$ is an iterated tensor product of factors of the form
$(D^-)^\ell(\si)$ and $(D^-)^\ell K$ for a component $K$ of $\ka$. The form of this
tensor product is required to be independent of $\si$ and it determines the
representation on which the map inducing $A_i$ is defined.

\begin{thm}\label{affine-invar}
Let $\bv$ and $\Bbb W$ be representations of $P$ and let $\Ph:\Ga(\Cal VA)\to\Ga(\Cal
WA)$ be a polynomial invariant operator constructed from $D^-$ in the sense
introduced above. Then for any Weyl structure $s$, the operator $\Ph^s:\Ga(\Cal
VM)\to\Ga(\Cal WM)$ is induced, under the isomorphism to the associated graded
bundles determined by $s$, by a polynomial invariant operator constructed from
$\nabla^s$, its curvature and its torsion in the same sense. Moreover,
this correspondence has the property that for two Weyl structures $s$ and $\tilde s$
and a point $x\in M$ such that $s(x)=\tilde s(x)$ one obtains
$\Ph^s(\si)(x)=\Ph^{\tilde s}(\si)(x)$ for any $\si\in\Ga(\Cal VM)$.
\end{thm}
\begin{proof}
  Theorem \ref{k-th-order} provides us with an expansion for $(D^-)^\ell(\si)$ for
  any $\ell$ involving $G_0$-equivariant maps acting on tensors. The discussion
  in Section \ref{involve-curv}
  gives us analogous expressions for $(D^-)^\ell(K)$ for any component $K$ of
  $\ka$. Since tensor products of $G_0$-equivariant maps are $G_0$-equivariant, we
  also get expressions for iterated tensor products of such terms. Applying
  $G_0$-equivariant maps to such tensor products just gives rise to compositions of
  $G_0$-equivariant maps, which are equivariant, too. The last statement has already
  been observed for any operator of the form $\Ph^s$ in the beginning of Section 2. 
\end{proof}

\section{The nearly invariant operators}

\subsection{The transformation rules}\label{transf}
The natural approach to invariants of parabolic geometries is via Weyl-structures,
i.e. to consider differential operators and differential invariants defined using the
Weyl connections $\nabla^s$, their torsion and their curvature, and request that they
are independent of the choice of $s$.  Thus, we should understand how the gradings,
the covariant derivatives and Rho tensors change under the change of the Weyl
structures.

Considering two Weyl structures as reductions $\bar s_1,\ \bar s_2:\cg_0\to\cg$,
clearly there must be a function $\Up:\cg_0\to \fp_+$ such that for all $u\in \cg_0$,
$g_0\in G_0$,
\begin{align*}
\bar s_2(u) &= \bar s_1(u)\cdot \operatorname{exp}\Up(u),
\\
\bar s_2(u\cdot g_0) &= \bar s_1(u)\cdot \operatorname{exp}\Up(u)\cdot g_0 = 
\bar s_1(u)\cdot g_0\cdot \operatorname{exp} (\operatorname{Ad}_{g_0^{-1}}
\Up(u)).
\end{align*}
Thus the function $\Up$ represents a one-form on $M$ and, fixing $s_1$, this is a
bijective correspondence between $\Up\in\Om^1(M)$ and all Weyl structures. This
reflects the affine bundle structure of $A$. Of course, these one-forms are also
represented as functions $\Up:\cg\to \fp_+$ with the right equivariancy property, and
then the functions $f=\operatorname{exp}\Up:\cg\to P/G_0\simeq P_+$ with the
equivariance property $f(u\cdot (g_0g_+))=g_+^{-1}g_0^{-1}f(u)g_0$ can be directly
seen as the corresponding sections $s:M\to A=\cg\x_P P/P_0$. See \cite[Proposition
  2.2]{CapMettler} for a detailed exposition.

The formulae for the transformations in terms of the
functions $\Up$ are explained in detail in \cite[Sections
5.1.6-9]{parabook}. They look pretty complicated and we shall not need them
explicitly here. Just in the simplest case of trivial filtrations (i.e.,
$|1|$-graded $\fg$) we obtain for vector fields $\xi$ on $M$ and sections
$\si$ of irreducible natural bundles (the hat indicates the transformed
objects)
\begin{align}
\hat\nabla_\xi \si &= \nabla_\xi \si - \{\Up,\xi\}\bullet \si
\label{change_nabla}\\
\hat\Rho(\xi) &= \Rho(\xi) + \nabla_\xi\Up
+\frac12\operatorname{ad}(\Up)^2(\xi) 
\label{change_rho}.\end{align}

Differentiating \eqref{change_nabla} again, derivatives of $\Up$ appear, and
the relatively nice transformation rule \eqref{change_rho} indicates a
chance to balance the formulae by adding `correction terms' 
including $\Rho$ in order to keep the transformation rules algebraic in
the parameters $\Up$. This is exactly what happens with the Rho-corrected
derivatives appearing in Proposition \ref{nabla_comparison}. 

About a hundred
years back this was the motivation for considering Schouten's Rho
tensor in conformal Riemannian geometry.

In the sequel we shall consider \emph{polynomial invariant operators}
$\Ps_s:\Gamma(\mathcal VM)\to \Gamma(\mathcal WM)$, constructed by a fixed universal
expression $\Ps$ from from $\nabla^s$, its curvature and its torsion in the way
discussed in Section \ref{affine-invar}. We say that $\Ps$ is a \emph{nearly
  invariant operator} if the transformation formula for $\Ps_s$ under the change of
the Weyl structure $s$ is tensorial in $\Up$.  Otherwise put, the operator $\Ps$ is
nearly invariant if and only if for Weyl structures $s_1,s_2$ and a point $x\in M$
such that $s_1(x)=s_2(x)$ we get $\Psi_{s_1}\si(x)=\Psi_{s_2}\si(x)$ for every
$\si\in\Gamma(\mathcal VM)$.

Theorem \ref{affine-invar} shows, that via $\Ph\mapsto\Ph^s$, any polynomial
affine invariant of the canonical connection $D^-$ on $A$ gives rise to a nearly
invariant operator on $M$.

\subsection{Derivatives of the Rho-tensor}\label{Rho-deriv}
The final key step towards understanding nearly invariant operators
  is analogous to the well known fact that in affine differential invariants one may
  use symmetrized iterated covariant derivatives rather then iterated covariant
  derivatives. We need an analogous symmetrization argument for iterated covariant
  derivatives of the Rho tensor associated to a Weyl structures. This will allow us
  to effectively use normal Weyl structures in the study of nearly invariant
  operators. 

Recall from Section \ref{normal_weyl} that a characteristic of normal Weyl structures
is the vanishing of certain symmetrizations of iterated covariant derivatives of the
Rho-tensor.  We will refer to these as \textit{form-symmetrized} iterated
derivatives.  Fixing a Weyl-structure $s$ with Rho-tensor $\Rho$ and Weyl connections
$\nabla$, we view $\Rho$ as a one-form on $M$ with values in $\gr(T^*M)$.  The $k$th
covariant derivative $(\nabla)^k\Rho$ is then a section of $\otimes^{k+1}T^*M\otimes
T^*M$ (or the associated graded of this bundle).  Now we define the $k$th
\textit{form symmetrized Weyl-derivative} $\Cal F^k(\Rho)\in
\Ga(\gr(S^{k+1}T^*M\otimes T^*M))$ as the symmetrization of $(\nabla^s)^k\Rho$ over
the first $k+1$ indices.  For convenience, we put $\Cal F^0(\Rho)=\Rho$.

\begin{thm}\label{thm-Rho-deriv}
  For any $k$, there is a universal expression for $(\nabla)^k\Rho-\Cal F^k(\Rho)$
  obtained via bundle maps induced by $G_0$-equivariant maps from tensor products
  whose factors are of the form $\Cal F^\ell(\Rho)$ and $(\nabla^s)^\ell K$ with
  $\ell<k$, where $K$ is a component of the pullback of the Cartan curvature along
  $s$.
\end{thm}
\begin{proof}
  Throughout this proof, we write $\nabla$ for $\nabla^s$.  We start by computing
  $(\nabla\Rho-\Cal F^1(\Rho))(\xi,\eta)$ for vector fields $\xi,\eta\in\frak X(M)$
  via
  $$
  \nabla_\xi\Rho(\eta) - \nabla_\eta\Rho(\xi)= \nabla_\xi(\Rho(\eta)) - \Rho(\nabla_\xi \eta) -
  \nabla_\eta(\Rho(\xi)) + \Rho(\nabla_\eta \xi)
  .$$
  The definition of the torsion of a Weyl connection reads as
  \begin{equation}\label{tors}
  T(\xi,\eta)=\nabla_\xi\eta-\nabla_\eta\xi-[\xi,\eta]+\{\xi,\eta\},
  \end{equation}
  so we can insert this to rewrite the terms in which a derivative goes into Rho. Moreover,
  from section \ref{curvature_comparison}, we know that
  $$
  Y^s(\xi,\eta)=\nabla_\xi(\Rho(\eta))-\nabla_\eta(\Rho(\xi))-\Rho([\xi,\eta])+
  \Rho(\{\xi,\eta\})+\{\Rho(\xi),\Rho(\eta)\}
  $$
  can be recoverd from the pullback of the Cartan curvature and a
    component of $\Rho$. Together with the above, we conclude that
  \begin{equation}\label{Rho-sym1}
  \nabla\Rho=\Cal F^1(\Rho)+\tfrac12 Y-\tfrac12 i_T\Rho-\tfrac12\{\Rho,\Rho\}.
\end{equation}

  This proves the theorem in the case $k=1$. To complete the proof we show by
  induction on $k$ that $\nabla \Cal F^k(\Rho)-\Cal F^{k+1}(\Rho)$ admits a universal
  expansion as claimed in theorem with $\ell\leq k$. Since, by induction, one Weyl
  derivative of any term allowed in step $k$ leads to a term allowed in step $k+1$,
  this recursively implies that claim of the theorem.

  The case $k=0$ has been sorted out above, so we assume that $k>0$ and our claim has
  been proved $\ell<k$.

Since $\Cal F^{k-1}(\Rho)\in\Ga(S^kT^*M\otimes T^*M)$ is symmetric in its first
$k$-indices, we can compute $\Cal F^k\Rho$ as $\tfrac{1}{k+1}$ times the sum over all
cyclic permutations in the first $k+1$ entries of $\nabla \Cal F^{k-1}(\Rho)\in
\Ga(T^*M\otimes S^kT^*M\otimes T^*M)$.
Using abstract index notation $\Cal F_{a_1\dots a_k b}$ for $\Cal F^{k-1}\Rho$ we can
write $\Cal F^k(P)_{a_0 \dots a_k b}$ as
  $$
\frac{1}{k+1}\sum_{i=0}^k\nabla_{a_i}\Cal F_{a_0\dots \widehat{a_i} \dots a_k b},
$$
where as usual the hat denotes omission and we observe we may arbitrarily permute the
first $k$ indices in $\Cal F$. Using symmetry of $\Cal F$ to rewrite
$\nabla_{a_0}\Cal F_{a_1 \dots a_k b}$, we conclude that we can rewrite $\nabla \Cal
F^k(\Rho)-\Cal F^{k+1}(\Rho)$ as
$$
\frac{1}{k+1}\sum_{i=1}^k\left(\nabla_{a_0}\Cal
F_{a_i a_1, \dots \widehat{a_i} \dots a_kb}-\nabla_{a_i}\Cal F_{a_0a_1
  \dots\widehat{a_i}\dots a_k b}\right). 
$$
If $k=1$, then the terms $\Cal F$ already are of the form $\nabla\Rho$. For $k>1$,
the induction hypothesis implies that, up to terms of the form we allow in our
expansions, we may replace the occurrences of $\Cal F$ by $\nabla\Cal
F^{k-2}(\Rho)$. But then each term in our sum becomes twice an alternation of a
double covariant derivative of $\Cal F^{k-2}$. (Recall that $\Cal F^0(\Rho)=\Rho$.)

But now for any tensor field $t$, the alternation of $\nabla^2t$ maps
$\xi,\eta\in\frak X(M)$ to
$$
\nabla_\xi\nabla_\eta t-\nabla_{\nabla_\xi\eta}t-\nabla_\eta\nabla_\xi
t+\nabla_{\nabla_\eta\xi}t. 
$$
Inserting from \eqref{tors} and rearranging terms, we see that this can be rewritten
as
\begin{equation}\label{alt-n2}
R(\xi,\eta)\bullet t - \nabla_{T(\xi,\eta)}t+\nabla_{\{\xi,\eta\}}t,
\end{equation}
where in the first term the bullet denotes the tensorial action of the curvature.
Now the last two terms in \eqref{alt-n2} are obtained from $G_0$-equivariant
operations acting on $T\otimes \nabla t$ and of $\nabla t$, respectively, and $T$
can be obtained from the pullback of the Cartan curvature and from
  $\Rho$. In our situation $t=\Cal F^{k-2}(\Rho)$, so by induction we can rewrite
$\nabla t$ as a sum of terms of the allowed forms, so we see that these two summands
altogether only produce terms of the allowed form.

For the first terms in \eqref{alt-n2}, we can rewrite the curvature $R$ of
$\nabla=\nabla^s$ according to section \ref{curvature_comparison} as $W-
  (\partial\Rho)_{\fg_0}$ where $W$ is a component of the pullback of the Cartan
curvature. But then the resulting terms can be obtained via $G_0$-equivariant maps
from $W\otimes t$ respectively from $\Rho\otimes t$ and since $t$ equals $\Cal
F^{k-2}(\Rho)$, this leads to allowed terms only.
\end{proof}

\subsection{Nearly invariant operators}\label{near-inv}
Let us collect the information on affine invariants of the Weyl connections we have
obtained so far. As discussed in Section \ref{transf}, such an operator comes from a
universal expression $\Psi$ which involves bundle maps induced from
$G_0$-equivariant linear maps and certain tensor products. If the operator acts on
$\Ga(\Cal VM)$ these tensor products contain factors that are iterated covariant
derivatives of either a section $\si\in\Ga(\Cal VM)$ or a component of the curvature
of the involved connection (for which we will insert all Weyl connections). Now as
discussed in Sections \ref{curvature_comparison} and \ref{involve-curv}, we can
universally decompose the curvature and torsion of the Weyl connections into
components of the pullback of the Cartan connection and into components of the
Rho-tensor. Since these decompositions are induced by $G_0$-equivariant maps, the
same applies to iterated derivatives of (components of) the torsion and the
curvature. Moreover, we can use Theorem \ref{thm-Rho-deriv} to rewrite iterated
derivatives of the Rho-tensor via components of the pullback of the Cartan curvature
and form-symmetrized iterated derivatives of the Rho-tensor and as before, this
extends to components of the Rho-tensor. The upshot of this is that we may assume
that all our factors in the tensor products are
\begin{itemize}
\item iterated covariant derivatives of $\si\in\Ga(\Cal VM)$
\item iterated covariant derivatives of components of the pullback of the Cartan
    curvature
\item components of form symmetrized iterated covariant derivatives of the Rho tensor.
\end{itemize}

Armed with this observation, we are ready to formulate and proof the main result of
the paper. 

\begin{thm}\label{main-thm}
The nearly invariant operators are exactly the universal expansions obtained via
$\Ph\mapsto\Ph^s$ from polynomial affine differential invariants of the natural
covariant derivative $D^-$ and its curvature and torsion on $A$.
\end{thm}
\begin{proof}
Theorem \ref{affine-invar} shows that affine differential invariants of $D^-$ and the
curvature and torsion of $D$ induce nearly invariant operators.

To prove the converse inclusion, let us take a universal expression $\Ps$ which gives
rise to a nearly invariant operator which is formed in the way discussed above. Now
let $\tilde\Ps^1$ be the universal expression obtained by removing from $\Ps$ all
terms in which the tensor product contains a component of a form symmetrized
covariant derivative of the Rho tensor. For each summand in this expression count the
total number of covariant derivatives in all factors showing up in the tensor
product. Let $k$ be the maximal number that occurs and let $\Ps_1$ be the sum of all
terms in which the total number of covariant derivatives equals $k$. Then there is an
obvious affine invariant $\Phi_1$ of $D^-$ on $A$ that corresponds to $\Psi_1$ (in
which we simply replace all covariant derivatives by $D^-$). Via $\Phi_1\mapsto
\Phi_1^s$, we obtain a nearly invariant operator defined on $\Ga(\Cal VM)$ and we
subtract this from $\Psi$.

On the one hand, the result by construction is a nearly invariant operator
$\tilde\Ps$. On the other hand, Theorem \ref{k-th-order} implies that in $\Phi_1^s$,
we obtain exactly the same terms involving a total number of $k$ covariant
derivatives as in $\Psi^1_s$ while for all other terms, the total number of covariant
derivatives that occur is strictly less than $k$. Otherwise put, the operator
$\tilde\Ps$ has the property that in any term that does not contain a component
of a form symmetrized iterated covariant derivative of the Rho tensor, there are less
than $k$ covariant derivatives in total.

This shows that we can iterate our procedure by applying the same construction to
$\tilde\Psi$. After finitely many steps, we arrive at an expression for which the
terms that do not involve any components of form symmetrized iterated covariant
derivatives of Rho also do not contain any covariant derivatives and hence evidently
are obtained from an affine invariant on $A$.

Hence we conclude that our original operator $\Psi$ can be written as the sum of an
operator of the form $\Ph\mapsto\Ph^s$ and a nearly invariant operator $\hat\Psi$ for
which \textit{any} term involves a component of a form symmetrized iterated covariant
derivative of Rho. But now we can easily conclude that the nearly invariant operator
$\hat\Psi$ vanishes identically. Indeed, take any Weyl structure $s$ and consider
$\Psi_s(\si)$ and a point $x\in M$. Then from Section \ref{normal_weyl} we know that
we can find a Weyl structure $\tilde s$ which is normal at $x$ and satisfies $\tilde
s(x)=s(x)$. But by normality, any term in $\hat\Psi_{\tilde s}$ vanishes at $x$ and
hence $0=\hat\Psi_{\tilde s}(\si)(x)$, and this equals $\hat\Psi_{s}(\si)(x)$ since
$\hat\Ps$ is nearly invariant. This completes the proof.  
\end{proof}

\subsection{Remark}
In our discussion, we have focused on the case of normal parabolic geometries, which
equivalently encode some underlying structure and hence provides applications to
invariants of that structure. The concept of Weyl structures as well as the basic
calculus we develop here are available for general Cartan geometries and basically
also the results we have proved above allow for extensions to this more general
setting. The basic difference is that the relation between the pullback $s^*\ka$ of
the Cartan curvature, the curvature of the Weyl connection $\nabla^s$ determined by
$s$ and the Rho-tensor $\Rho^s$ becomes more complicated. So to extend the results,
one would have to modify the definitions of affine invariants that are used in
Theorems \ref{affine-invar} and \ref{main-thm} appropriately. We do not go into
details on this here.

\end{document}